# Cartan equivalence problem for third order differential operators

M. Nadjafikhah [*], R. Bakhshandeh-Chamazkoti,

*School of Mathematics, Iran University of Science and Technology, Narmak, Tehran 1684613114, Iran.*

**Abstract**

This article is dedicated to solve the equivalence problem for two third-order differential operators on the line under general fiber-preserving transformation using the Cartan method of equivalence. We will do two versions of the equivalence problems: first via the direct equivalence problem, second equivalence problem is to determine conditions on two differential operators such that there exists a fiber-preserving transformations mapping one to the other according to gauge equivalence.

*Key words:* differential operator, Cartan equivalence, gauge equivalence, invariant, pseudogroup, Lie algebra.

## 1 Introduction

The classification of linear differential equations is a special case of the general problem of classifying differential operators, which has a variety of important applications, including quantum mechanics and the projective geometry of curves [1]. In this attempt we shall solve the method of local equivalence problem by three versions of the equivalence problem for the class of linear third order operators on the line. For simplicity, we shall only deal with the local equivalence problem for scalar differential operators in a single independent variable, although these problems are important for matrix–valued and partial differential operators as well. The general equivalence problem is to recognize when two geometrical objects are mapped on each other by a certain class of diffeomorphisms. E. Cartan developed the general equivalence problem and provided a systematic procedure for determining the necessary and sufficient condition [2,3]. In Cartan's approach, the conditions of equivalence of two objects must be reformulated in terms of differential forms. We associate a collection of one-forms to an object under investigation in the original coordinates; the corresponding object in the new coordinates will have its own collection of one-forms. Once an equivalence problem has been reformulated in the proper Cartan form, in terms of a coframe $\omega$ on the $m$-dimensional base manifold $M$, along with a structure group $G \subset \mathrm{GL}(m)$, we can apply the Cartan equivalence method. The goal is to normalize the structure group valued coefficients in a suitably invariant manner, and this is accomplished through the determination of a sufficient number of invariant combinations thereof [1].

Recently, there are some works on solving equivalence problem on third order ODEs, [4,5]. The problems here are related to the more general equivalence problem for third order ordinary differential equations which E. Cartan studied under point transformations [6], and S. S. Chern turned his attention to the problem under contact transformations [7] and Hajime Sato et all [8], but are specialized by linearity. Niky Kamran and Peter J. Olver have been solved equivalence problem for second order differential operator with two versions of the equivalence problem [9].

---

[*] Tel. +9821-73913426. Fax +9821-77240472.

   *Email addresses:* m_nadjafikhah@iust.ac.ir (M. Nadjafikhah), r_bakhshandeh@iust.ac.ir (R. Bakhshandeh-Chamazkoti).

## 2 Equivalence of third order differential operators

Our starting point is a third order differential operator

$$\mathscr{D} = \sum_{i=0}^{3} f_i(x)\, D^i \tag{2.1}$$

where $f_0, f_1, f_2$ and $f_3$ are analytic functions of the real variable $x$, $D^i = d/dx^i$, and $D^0 = \mathrm{Id}$ be the identity operator. With applying $\mathscr{D}$ on a scalar-valued function $u(x)$, we obtain the following expression

$$\mathscr{D}[u] = \sum_{i=0}^{3} f_i(x)\, D^i u. \tag{2.2}$$

We discuss the problem of equivalence under general fiber-preserving transformations which are linear in the dependent variable

$$\bar{x} = \xi(x), \qquad \bar{u} = \varphi(x)\, u, \tag{2.3}$$

where $\varphi(x) \neq 0$. The total derivative operators are related by the chain rule formula

$$\bar{D} = \frac{d}{d\bar{x}} = \frac{1}{\xi'(x)} \frac{d}{dx} = \frac{1}{\xi'(x)} D. \tag{2.4}$$

We first consider the *direct equivalence problem*, which identifies the two linear differential functions

$$\sum_{i=0}^{3} f_i(x)\, D^i u = \mathscr{D}[u] = \bar{\mathscr{D}}[\bar{u}] = \sum_{i=0}^{3} \bar{f}_i(\bar{x})\, \bar{D}^i \bar{u}. \tag{2.5}$$

under the change of variables (2.3). This induces the transformation rule

$$\bar{\mathscr{D}} = \mathscr{D} \cdot \frac{1}{\varphi(x)} \qquad \text{when} \qquad \bar{x} = \xi(x), \tag{2.6}$$

on the differential operators themselves, and we try to find explicit conditions on the coefficients of the two differential operators that guarantee that they satisfy (2.5) for some change of variables of the form (2.3).

The transformation rule (2.6) doesn't preserve either the eigenvalue problem $\mathscr{D}[u] = \lambda u$ or the Schrödinger equation $i u_t = \mathscr{D}[u]$, since we are missing a factor of $\varphi(x)$. To rectify this problem, we need to multiply by $\varphi(x)$ and use the *gauge equivalence* with the following transformation rule

$$\bar{\mathscr{D}} = \varphi(x) \cdot \mathscr{D} \cdot \frac{1}{\varphi(x)} \qquad \text{when} \qquad \bar{x} = \xi(x). \tag{2.7}$$

In quantum mechanics, equivalence plays an important role since it preserves the solution set to the associated Schrödinger equation, or its stationary counterpart, the eigenvalue problem.

The appropriate space to work in will be the third jet space $\mathrm{J}^3$, which has local coordinates $\Upsilon = \{(x, u, p, q, r) \in \mathrm{J}^3 : p = u_x, q = u_{xx}, r = u_{xxx}\}$, and our goal is to construct an appropriate coframe on $\mathrm{J}^3$, which will encode the relevant transformation rules for our problem. Note first that a point transformation will be in the desired linear form (2.3) if and only if, for some pair of functions $\alpha, \beta$, one-form equations

$$d\bar{x} = \alpha\, dx, \tag{2.8}$$

$$\frac{d\bar{u}}{\bar{u}} = \frac{du}{u} + \beta\, dx. \tag{2.9}$$

hold on the subset of J³ where $u \neq 0$. Indeed, the (2.8) equation implies that $\bar{x} = \xi(x)$, with $\alpha = \xi_x$, while the (2.9) necessarily requires $\bar{u} = \varphi(x)u$, with $\beta = \varphi_x/\varphi$.

In order that the derivative variables $p, q$ and $r$ transform correctly, we need to preserve the contact ideal $\mathcal{I}$ on J³, which is

$$\mathcal{I} = \langle du - p\,dx, dp - q\,dx, dq - r\,dx \rangle. \tag{2.10}$$

Generally, a diffeomorphism $\Phi : \text{J}^3 \to \text{J}^3$ determines a contact transformation if and only if

$$d\bar{u} - \bar{p}\,d\bar{x} = a_1(du - p\,dx), \tag{2.11}$$
$$d\bar{p} - \bar{q}\,d\bar{x} = a_2(du - p\,dx) + a_3(dp - q\,dx), \tag{2.12}$$
$$d\bar{q} - \bar{r}\,d\bar{x} = a_4(du - p\,dx) + a_5(dp - q\,dx) + a_6(dq - r\,dx), \tag{2.13}$$

where $a_{ij}$ are functions on J³. The combination of the first contact condition (2.11) with the linearity conditions (2.8) and (2.9) constitutes part of an overdetermined equivalence problem. We put $\beta = -p/u$ but say before that $\beta = \phi_x/\phi$. It is better to simply to say that (2.10, 2.9) is equivalent to condition that form $\omega_2$ doesn't change under point preserving transformation

$$\frac{d\bar{u} - \bar{p}\,d\bar{x}}{\bar{u}} = \frac{du - p\,dx}{u}, \tag{2.14}$$

which is invariant, and (2.14) can replace both (2.9) and (2.11). Therefore, we choose four elements of our coframe the one-forms

$$\omega^1 = dx, \quad \omega^2 = \frac{du - p\,dx}{u}, \quad \omega^3 = dp - q\,dx, \quad \omega^4 = dq - r\,dx. \tag{2.15}$$

which are defined on the third jet space J³ locally parameterized by $(x, u, p, q, r)$, with the transformation rules

$$\bar{\omega}^1 = a_1\omega^1, \quad \bar{\omega}^2 = \omega^2, \quad \bar{\omega}^3 = a_2\omega^2 + a_3\omega^3, \quad \bar{\omega}^4 = a_4\omega^2 + a_5\omega^3 + a_6\omega^4. \tag{2.16}$$

Viewing to (2.5), the function $I(x, u, p, q, r) = \mathscr{D}[u] = f_3(x)r + f_2(x)q + f_1(x)p + f_0(x)u$ is an invariant for the problem, and thus its differential

$$\omega^5 = dI = f_3 dr + f_2 dq + f_1 dp + f_0 du + (f_3'r + f_2'q + f_1'p + f_0'u)dx, \tag{2.17}$$

is an invariant one-form, thus we take it as final element of our coframe.

In the second problem (2.7), for the extra factor of $\varphi$, the invariant is

$$I(x, u, p, q, r) = \frac{\mathscr{D}[u]}{u} = \frac{f_3(x)r + f_2(x)q + f_1(x)p}{u} + f_0(x). \tag{2.18}$$

Thus, we take

$$\omega^5 = dI = \frac{f_3}{u}dr + \frac{f_2}{u}dq + \frac{f_1}{u}dp - \frac{f_3 r + f_2 q + f_1 p}{u^2}du + \left\{\frac{f_3'r + f_2'q + f_1'p}{u} + f_0'\right\}dx, \tag{2.19}$$

as a final element of coframe for the second equivalence problem (2.7). In the both cases, the set of one-forms $\{\omega^1, \omega^2, \omega^3, \omega^4, \omega^5\}$ is a coframe on the subset

$$\Omega^* = \left\{ (x, u, p, q, r) \in \text{J}^3 \,\middle|\, u \neq 0 \text{ and } f_3(x) \neq 0 \right\}. \tag{2.20}$$

We restrict our attention to a connected component $\Omega \subset \Omega^*$ of the subset (2.20) that the signs of $f_0(x)$ and $u$ are fixed. In the both first and second problems, since $\omega^5 = dI$ is a closed invariant one-form, the last coframe elements agree up to

$$\bar{\omega}^5 = \omega^5. \tag{2.21}$$

In view of (2.16) and (2.21), the structure group associated with the equivalence problem is the six–dimensional matrix group

$$G = \left\{ \begin{pmatrix} a_1 & 0 & 0 & 0 & 0 \\ 0 & 1 & 0 & 0 & 0 \\ 0 & a_2 & a_3 & 0 & 0 \\ 0 & a_4 & a_5 & a_6 & 0 \\ 0 & 0 & 0 & 0 & 1 \end{pmatrix} : a_i \in \mathbb{R};\ i = 1, \cdots, 6,\quad a_1 a_3 a_6 \neq 0 \right\} \tag{2.22}$$

which happens to be the same for both equivalence problems, direct equivalence and gauge equivalence methods, even though the two coframes are different.

**Proposition 2.1** *Suppose $\mathscr{D}$ and $\bar{\mathscr{D}}$ be third-order differential operators. Let $\{\omega^1, \omega^2, \omega^3, \omega^4, \omega^5\}$, and $\{\bar{\omega}^1, \bar{\omega}^2, \bar{\omega}^3, \bar{\omega}^4, \bar{\omega}^5\}$ be the corresponding coframes, on open subsets $\Omega$ and $\bar{\Omega}$ of the third jet space, given by (2.15) and (2.17) or (2.19), the choice of $\omega^5$ and $\bar{\omega}^5$ depending on the equivalence problem under consideration. The differential operators are equivalent under the pseudogroup (2.3) according to the respective transformation rule (2.6) and (2.7) if and only if there is a diffeomorphism $\Phi$ that satisfies*

$$\Phi^*(\bar{\omega}_i) = \sum_{j=1}^{5} g_{ij}\, \omega_j, \tag{2.23}$$

*for $i = 1, \cdots, 5$, where $g = (g_{ij})$ is a $G$–valued function on $\mathrm{J}^3$, and $\Phi^*$ denotes the pull-back map on differential forms.*

In order to apply Cartan's reduction algorithm for direct quivalence and gauge equivalence problems so as to prescribe invariant normalizations of the six group parameters $a_1, a_2, \cdots, a_6$, we must lift coframes to the space $\mathrm{J}^3 \times G$. The *lifted coframe* has the form

$$\begin{aligned}
\theta^1 &= a_1 \omega^1, \\
\theta^2 &= \omega^2, \\
\theta^3 &= a_2 \omega^2 + a_3 \omega^3, \\
\theta^4 &= a_4 \omega^2 + a_5 \omega^3 + a_6 \omega^4, \\
\theta^5 &= \omega^5.
\end{aligned} \tag{2.24}$$

By normalizing the torsion components, we have managed to eliminate all of the group parameters. This has had the effect of reducing the structure group to the identity, and reducing the lifted invariant coframe to an invariant coframe on the base space $\mathrm{J}^3$, known as an $\{e\}$–*structure* or *local parallelism*.

## 3 The direct equivalence problem

First, we take the initial four one-forms (2.15) and (2.17) as our final coframe constituent. The next step is to calculate the differentials of lifted coframe elements (2.24). An explicit computation leads to the structure equations

$$\begin{aligned}
d\theta^1 &= \alpha^1 \wedge \theta^1, \\
d\theta^2 &= T^2_{12} \theta^1 \wedge \theta^2 + T^2_{13} \theta^1 \wedge \theta^3, \\
d\theta^3 &= \alpha^2 \wedge \theta^2 + \alpha^3 \wedge \theta^3 + T^3_{12} \theta^1 \wedge \theta^2 + T^3_{13} \theta^1 \wedge \theta^3 + T^3_{14} \theta^1 \wedge \theta^4, \\
d\theta^4 &= \alpha^4 \wedge \theta^2 + \alpha^5 \wedge \theta^3 + \alpha^6 \wedge \theta^4 + T^4_{12} \theta^1 \wedge \theta^2 + T^4_{13} \theta^1 \wedge \theta^3 + T^4_{14} \theta^1 \wedge \theta^4 + T^4_{15} \theta^1 \wedge \theta^5, \\
d\theta^5 &= 0,
\end{aligned} \tag{3.1}$$

with

$$\alpha^1 = \frac{da_1}{a_1}, \qquad \alpha^2 = \frac{a_3 da_2 - a_2 da_3}{a_3}, \qquad \alpha^3 = \frac{da_3}{a_3},$$
$$\alpha^4 = \frac{a_3 a_6 da_4 + (a_2 a_5 - a_3 a_4) da_6 - a_2 a_6 da_5}{a_3 a_6}, \qquad \alpha^5 = \frac{a_3 da_5 - a_5 da_6}{a_3 a_6}, \qquad \alpha^6 = \frac{da_6}{a_6},$$

forming a basis for the right-invariant *Maurer-Cartan forms* on the Lie group $G$. The torsion coefficients in the structure equations (3.1) are explicitly given by

$$\begin{aligned}
T^2_{12} &= -\frac{a_2 + a_3 p}{a_1 a_3 u}, \quad T^2_{13} = \frac{1}{a_1 a_3 u}, \quad T^3_{12} = -\frac{a_2 a_6 a_3 p + a_2^2 a_6 - a_2 a_5 a_3 u + a_3^2 a_4 u}{a_1 a_3 a_6 u}, \\
T^3_{13} &= \frac{a_2 a_6 - a_3 a_5 u}{a_1 a_3 a_6 u}, \quad T^3_{14} = \frac{a_3}{a_1 a_6}, \\
T^4_{12} &= -\frac{(a_3 a_4 a_6 p + a_2 a_4 a_6 + a_4 a_3 a_5 - a_2 a_5^2) f_3 + (a_2 a_5 a_6 - a_4 a_3 a_6) f_2 - a_2 a_6^2 f_1 + a_6^2 a_3 u^2 f_0}{a_1 a_3 a_6 f_3 u}, \\
T^4_{13} &= \frac{(a_4 a_6 - a_5^2 u) f_3 + a_5 a_6 u f_2 - a_6^2 u f_1}{a_1 a_3 a_6 f_3 u}, \quad T^4_{14} = \frac{a_5 f_3 - a_6 f_2 - a_6 a_9}{a_1 a_6 f_3}, \quad T^4_{15} = \frac{a_6}{a_1 f_3}.
\end{aligned} \qquad (3.2)$$

In the absorption part of Cartan's process, we replace each Maurer-Cartan form in the structure equations by general linear combination of coframe elements, so $\alpha^\kappa \mapsto \sum_{j=1}^5 z_j^\kappa \theta^j$, where the coefficients $z_j^\kappa$ are allowed to depend on both the base variables $x, u, p, q, r$ and the group parameters $a_1, a_2, \ldots, a_6$. The resulting two-forms are

$$\begin{aligned}
\Theta^1 &= -z_2^1 \theta^1 \wedge \theta^2 - z_3^1 \theta^1 \wedge \theta^3 - z_4^1 \theta^1 \wedge \theta^4 - z_5^1 \theta^1 \wedge \theta^5, \\
\Theta^2 &= T^2_{12} \theta^1 \wedge \theta^2 + T^2_{13} \theta^1 \wedge \theta^3, \\
\Theta^3 &= (z_1^2 + T^3_{12}) \theta^1 \wedge \theta^2 + (z_1^3 + T^3_{13}) \theta^1 \wedge \theta^3 + T^3_{14} \theta^1 \wedge \theta^4 + (z_2^3 - z_3^2) \theta^2 \wedge \theta^3 \\
&\quad - z_4^2 \theta^2 \wedge \theta^4 - z_4^3 \theta^3 \wedge \theta^4 - z_5^2 \theta^2 \wedge \theta^5 - z_5^3 \theta^3 \wedge \theta^5, \\
\Theta^4 &= (z_1^4 + T^4_{12}) \theta^1 \wedge \theta^2 + (z_1^5 + T^4_{13}) \theta^1 \wedge \theta^3 + (z_1^6 + T^4_{14}) \theta^1 \wedge \theta^4 + (z_2^5 - z_3^4) \theta^2 \wedge \theta^3 \\
&\quad + (z_2^6 - z_4^4) \theta^2 \wedge \theta^4 + (z_3^6 - z_4^5) \theta^3 \wedge \theta^4 + T^4_{15} \theta^1 \wedge \theta^5 - z_5^4 \theta^2 \wedge \theta^5 - z_5^5 \theta^3 \wedge \theta^5 - z_5^6 \theta^4 \wedge \theta^5, \\
\Theta^5 &= 0.
\end{aligned}$$

Some coefficients of $\theta^j \wedge \theta^k$ in each $\Theta^i$ which is independent of the parameters $z_j^\kappa$ are invariants of the problem, so one can normalize to reduce the structure group. In above, the essential torsion components are $T^2_{12}, T^2_{13}, T^3_{14}, T^4_{15}$, as given in (3.2) which is possible to absorb all the torsion components except them. By direct inspection of the structure equations (3.1), we deduce that any torsion components in $d\theta^2$ are essential because there are no Maurer-Cartan forms in it and since the Maurer-Cartan forms in $d\theta^3$ multiply either $\theta^2$ or $\theta^3$, and $d\theta^4$ multiply $\theta^2, \theta^3$ and $\theta^4$ they can never produce a multiple of the two-form $\theta^1 \wedge \theta^4$ and $\theta^1 \wedge \theta^5$ upon replacement respectively.

Since the essential torsion coefficients all depend on the group parameters, thus the next step in the process is to normalize them to as simple a form as possible. We first normalize $T^2_{12} = 0$ by setting $a_2 = -a_3 p$, thereby eliminating the group parameter $a_2$. Second, we can normalize $T^3_{14} = T^4_{15} = 1$ by setting $a_3 = a_1 a_6$, $a_6 = a_1 f_3$. With these four normalizations, the fourth essential torsion coefficient becomes $T^2_{13} = 1/(a_1^3 f_3 u)$. By assumption $f_3(x) u \neq 0$ and because of using real-valued functions thus $T^2_{13} = 1$. Therefore, we normalize $a_1 = (f_3 u)^{-1/3}$. The group parameter normalizations are

$$a_1 = \frac{1}{\sqrt[3]{f_3 u}}, \quad a_2 = -\sqrt[3]{\frac{f_3}{u^2}}\, p, \quad a_3 = \sqrt[3]{\frac{f_3}{u^2}}, \quad a_6 = \sqrt[3]{\frac{f_3}{u}}. \qquad (3.3)$$

Now substituting normalizations (3.3) in the lifted coframe (2.24), in the second loop through the equivalence procedure, we calculate the differentials of new invariant coframe and so the revised structure equations are

$$
\begin{aligned}
d\theta^1 &= \frac{1}{3}\theta^1 \wedge \theta^2, \\
d\theta^2 &= \theta^1 \wedge \theta^3, \\
d\theta^3 &= T^3_{12}\theta^1 \wedge \theta^2 + T^3_{13}\theta^1 \wedge \theta^3 + \theta^1 \wedge \theta^4 + \frac{1}{3}\theta^2 \wedge \theta^3, \\
d\theta^4 &= \alpha^1 \wedge \theta^2 + \alpha^2 \wedge \theta^3 + T^4_{12}\theta^1 \wedge \theta^2 + T^4_{13}\theta^1 \wedge \theta^3 + T^4_{14}\theta^1 \wedge \theta^4 + T^4_{23}\theta^2 \wedge \theta^3 \\
&\quad -\frac{1}{3}\theta^2 \wedge \theta^4 + \theta^1 \wedge \theta^5, \\
d\theta^5 &= 0,
\end{aligned}
\tag{3.4}
$$

where $\alpha^1$ and $\alpha^2$ are the Maurer-Cartan forms on the structure group $G$. The essential torsion coefficients are

$$
T^3_{12} = -a_4 - \sqrt[3]{\frac{f_3^2}{u}}\, q, \quad T^3_{13} = \frac{uf_3' - 5pf_3}{3\sqrt[3]{(f_3u)^2}} - a_5, \quad T^4_{14} = a_5 + \frac{2f_3'u - 3f_2u - f_3p}{3\sqrt[3]{(f_3u)^2}}.
\tag{3.5}
$$

Since in (3.4) the other torsion coefficients can be absorbed by the Maurer-Cartan forms, we just normalize the essential torsion coefficients (3.5) and we have found the following parameters:

$$
a_4 = -\sqrt[3]{\frac{f_3^2}{u}}, \qquad a_5 = \frac{f_3'u - 5f_3p}{3\sqrt[3]{(f_3u)^2}}.
\tag{3.6}
$$

The normalizations (3.6) have the effect of reducing the original structure group $G$ to a one-parameter subgroup and we have finally normalized all the group parameters. Inserting their prescribed values (3.3) and (3.6) into (2.24), the invariant coframe is now given by

$$
\begin{aligned}
\theta^1 &= \frac{dx}{\sqrt[3]{f_3u}}, \\
\theta^2 &= \frac{du - p\,dx}{u}, \\
\theta^3 &= \sqrt[3]{\frac{f_3}{u^2}}\left[(dp - q\,dx) - \frac{p}{u}(du - p\,dx)\right], \\
\theta^4 &= -\sqrt[3]{\frac{f_3^2}{u}}\Big(\frac{qdu - pq\,dx}{u}\Big) + \frac{f_3'u - 5f_3p}{3\sqrt[3]{(f_3u)^2}}(dp - q\,dx) + \sqrt[3]{\frac{f_3}{u}}(dq - r\,dx), \\
\theta^5 &= f_3dr + f_2dq + f_1dp + f_0du + (f_3'r + f_2'q + f_1'p + f_0'u)dx.
\end{aligned}
\tag{3.7}
$$

The final structure equations are then found to be

$$
\begin{aligned}
d\theta^1 &= \frac{1}{3}\theta^1 \wedge \theta^2, \\
d\theta^2 &= \theta^1 \wedge \theta^3, \\
d\theta^3 &= \theta^1 \wedge \theta^4 + \frac{1}{3}\theta^2 \wedge \theta^3, \\
d\theta^4 &= -I\theta^1 \wedge \theta^2 + \frac{1}{9}I_2\theta^1 \wedge \theta^3 + I_1\theta^1 \wedge \theta^4 + \theta^1 \wedge \theta^5 + \frac{2}{3}\theta^2 \wedge \theta^4, \\
d\theta^5 &= 0,
\end{aligned}
\tag{3.8}
$$

where the functions

$$\begin{aligned}
I &= f_3 r + f_2 q + f_1 p + f_0 u, \\
I_1 &= \frac{(f_3' - f_2)u - 2 f_3 p}{3 \sqrt[3]{(f_3 u)^2}}, \\
I_2 &= \frac{1}{f_3 u \sqrt[3]{f_3 u}} \Big[ (3 f_3' f_2 + 3 f_3 f_3'' - 4 f_3'^2 - 9 f_1 f_3) u^2 + 5 f_3^2 p^2 - 24 f_3^2 u q + (7 f_3' f_3 - 15 f_3 f_2) u p \Big],
\end{aligned} \quad (3.9)$$

are the fundamental structure invariants of the problem.

In a local coordinates $\Upsilon$ on $J^3$, the coframe can be written in terms of the coordinate coframe, so that $\theta^i = \sum_j a_j^i(x) dx^j$, where $A = (a_j^i(x))$ is a nonsingular $m \times m$ matrix of functions. The differential of a function can be re-write in the coframe-adapted form

$$dF = \sum_{j=1}^{5} \frac{\partial F}{\partial \theta^j} \theta^j, \quad (3.10)$$

which one will refer to the resulting coefficients $\partial F / \partial \theta^j = \sum_i b_j^i(x) \partial F / \partial x^i$, as the *coframe derivatives* of the function $F$, where $B = (b_j^i(x)) = A^{-1}$. Comparing (3.10) and the formulae (3.7) for the invariant coframe, we find that coframe derivatives of a function $F$ are given explicitly by

$$\begin{aligned}
\frac{\partial F}{\partial \theta^1} &= \sqrt[3]{f_3 u} \, \widehat{D}_x F, \\
\frac{\partial F}{\partial \theta^2} &= u \frac{\partial F}{\partial u} + p \frac{\partial F}{\partial p} + q \frac{\partial F}{\partial q} - \Big( \frac{f_2 q + f_1 p + f_0 u}{f_3} \Big) \frac{\partial F}{\partial r}, \\
\frac{\partial F}{\partial \theta^3} &= \frac{u}{\sqrt[3]{f_3 u}} \frac{\partial F}{\partial p} + \Big( \frac{5 p f_3 - f_3' u}{3 f_3 \sqrt[3]{f_3 u}} \Big) \frac{\partial F}{\partial q} - \Big( \frac{5 f_2 f_3 p - f_2 f_3' u + 3 f_1 f_3 u}{3 f_3^2 \sqrt[3]{f_3 u}} \Big) \frac{\partial F}{\partial r}, \\
\frac{\partial F}{\partial \theta^4} &= \frac{u}{\sqrt[3]{(f_3 u)^2}} \frac{\partial F}{\partial q} - \frac{f_2 u}{f_3 \sqrt[3]{(f_3 u)^2}} \frac{\partial F}{\partial r}, \\
\frac{\partial F}{\partial \theta^5} &= \frac{1}{f_3} \frac{\partial F}{\partial r},
\end{aligned} \quad (3.11)$$

where

$$\widehat{D}_x = \frac{\partial}{\partial x} + p \frac{\partial}{\partial u} + q \frac{\partial}{\partial p} + r \frac{\partial}{\partial q} + R \frac{\partial}{\partial r}, \quad R = -\frac{f_2 r + f_1 q + f_0 p + f_3' r + f_2' q + f_1' p + f_0' u}{f_3}.$$

The Jacobi identities for the coframe derivatives are found by reapplying the exterior derivative to the structure equations (3.8). An easy calculation shows that $d^2 \theta^1 = d^2 \theta^2 = d^2 \theta^3 = 0$ automatically, while the identities $d^2 \theta^4 = 0$ imply the following syzygy among our fundamental invariants and the derived invariants:

$$\frac{\partial I_1}{\partial \theta^2} = -I_1, \qquad \frac{\partial I_1}{\partial \theta^3} = 2, \qquad \frac{\partial I_2}{\partial \theta^2} = \frac{2}{3} I_2, \qquad \frac{\partial I_2}{\partial \theta^3} = 15 I_1, \qquad \frac{\partial I_2}{\partial \theta^4} = -24,$$

and also we have $dI = \theta^5$, meaning $\partial I / \partial \theta^5 = 1$ and $\partial I / \partial \theta^i = 0$ for $i = 1, 2, 3, 4$.

## 4 The gauge equivalence problem

We determine the solution to the problem of gauge equivalence of third differential operators by a similar computation of previous section. The Cartan formulation of this problem will use the same initial four one-forms (2.15), but now the final coframe element is (2.19). In the first loop through the second equivalence problem procedure, the structure group (2.22) is exactly the preliminary structure equations (3.1), and so the equivalence method has the same intrinsic structure by the unabsorbable torsion coefficients

$$T^2_{12} = -\frac{a_2 + a_3 p}{a_1 a_3 u}, \quad T^2_{13} = \frac{1}{a_1 a_3 u}, \quad T^3_{14} = \frac{a_3}{a_1 a_6}, \quad T^4_{15} = \frac{a_6 u}{a_1 f_3}. \tag{4.1}$$

It is possible to normalize the group parameters by setting

$$a_1 = \frac{1}{\sqrt[3]{f_3}}, \quad a_2 = -\sqrt[3]{\frac{f_3}{u^3}}\, p, \quad a_3 = \sqrt[3]{\frac{f_3}{u^3}}, \quad a_6 = \sqrt[3]{\frac{f_3^2}{u^3}}. \tag{4.2}$$

In the second loop of present equivalence problem, we substitute the normalizations (4.2) in the lifted coframe (2.23) and calculate the differentials of new invariant coframe to obtain revised structure equations. Now, we normalize the essential torsion components (4.1) by the remaind parameters

$$a_4 = -\sqrt[3]{\frac{f_3^2}{u^3}}\, q, \quad a_5 = \frac{f_3' u - 6 f_3 p}{3 u \sqrt[3]{f_3^2}}. \tag{4.3}$$

Thus the final invariant coframe is now given by

$$\begin{aligned}
\theta^1 &= \frac{dx}{\sqrt[3]{f_3}}, \\
\theta^2 &= \frac{du - p\, dx}{u}, \\
\theta^3 &= \frac{\sqrt[3]{f_3}}{u^2}\Big[(p^2 - qu)\, dx - p\, du + u\, dp\Big], \\
\theta^4 &= -\frac{1}{3\sqrt[3]{f_3}\, u^3}\Big[(3 f_3 u^2 r + f_3' u^2 q - f_3' u p^2 - 9 f_3 u p q + 6 f_3 p^3)\, dx + \\
&\qquad (f_3' u p + 3 f_3 u q - 6 f_3 p^2)\, du + (6 f_3 p - f_3' u) u\, dp - 3 f_3 u^2\, dq\Big], \\
\theta^5 &= \frac{f_3' r + f_2' q + f_1 p + f_0' u}{u}\, dx - \frac{f_3 r + f_2 q + f_1 p}{u^2}\, du + \frac{f_1}{u}\, dp + \frac{f_2}{u}\, dq + \frac{f_3}{u}\, dr.
\end{aligned} \tag{4.4}$$

The final structure equations take a following form

$$\begin{aligned}
d\theta^1 &= 0, \\
d\theta^2 &= \theta^1 \wedge \theta^3, \\
d\theta^3 &= \theta^1 \wedge \theta^4, \\
d\theta^4 &= I_1 \theta^1 \wedge \theta^3 + I_2 \theta^1 \wedge \theta^4 + \theta^1 \wedge \theta^5, \\
d\theta^5 &= 0,
\end{aligned} \tag{4.5}$$

where

$$\begin{aligned}
I_1 &= \frac{(f_3 f_3'' - 3 f_1 f_3 + f_2 f_3' - \frac{4}{3} f_3'^2) u + 3(f_3 f_3' - 2 f_2 f_3) p - 9 f_3^2 q}{3 f_3 \sqrt[3]{f_3}\, u}, \\
I_2 &= \frac{f_3' u - 3 f_3 p - f_2 u}{3 \sqrt[3]{f_3^2}\, u},
\end{aligned} \tag{4.6}$$

are the fundamental invariants of the problem. Note that the original invariant $I$, given in (2.18), does not appear among the structure functions of the adapted coframe. Nor can it appear among the derived invariants, since the

coframe derivatives are

$$\begin{aligned}
\frac{\partial F}{\partial \theta^1} &= \sqrt[3]{f_3}\, \widehat{D}_x F, \\
\frac{\partial F}{\partial \theta^2} &= u\frac{\partial F}{\partial u} + p\frac{\partial F}{\partial p} + q\frac{\partial F}{\partial q} + r\frac{\partial F}{\partial r}, \\
\frac{\partial F}{\partial \theta^3} &= \frac{u}{\sqrt[3]{f_3}}\frac{\partial F}{\partial p} + \Big(\frac{6f_3 p - f_3' u}{3 f_3 \sqrt[3]{f_3}}\Big)\frac{\partial F}{\partial q} + \Big(\frac{(f_2 f_3' - 3 f_1 f_3)u - 6 f_2 f_3 p}{3 f_3^2 \sqrt[3]{f_3}}\Big)\frac{\partial F}{\partial r}, \\
\frac{\partial F}{\partial \theta^4} &= \frac{u}{\sqrt[3]{f_3^2}}\frac{\partial F}{\partial q} - \frac{f_2 u}{f_3 \sqrt[3]{f_3^2}}\frac{\partial F}{\partial r}, \\
\frac{\partial F}{\partial \theta^5} &= \frac{u}{f_3}\frac{\partial F}{\partial r}.
\end{aligned} \qquad (4.7)$$

where

$$\widehat{D}_x = \frac{\partial}{\partial x} + p\frac{\partial}{\partial u} + q\frac{\partial}{\partial p} + r\frac{\partial}{\partial q} + R\frac{\partial}{\partial r}, \qquad (4.8)$$

$$R = -\frac{f_2 r u + f_1 q u - f_1 p^2 - f_2 p q - f_3 p r + f_3' r u + f_2' q u + f_1' p + f_0' u^2}{f_3 u}. \qquad (4.9)$$

The identity $d^2\theta^4 = 0$ leads to following syzygy among fundamental invariants and the derived invariants:

$$\frac{\partial I_1}{\partial \theta^4} = -\frac{\partial I_2}{\partial \theta^3} = -3, \qquad \frac{\partial I_1}{\partial \theta^3} = -2I_2. \qquad (4.10)$$